\documentclass[journal]{IEEEtran}
\usepackage{graphicx}
\usepackage{amsmath,amsfonts}
\usepackage{algorithmic}

\usepackage{amssymb}
\usepackage{textcomp,gensymb}
\usepackage{array}
\usepackage[caption=false,font=normalsize,labelfont=sf,textfont=sf]{subfig}
\usepackage{textcomp}
\usepackage{stfloats}
\usepackage{url}
\usepackage{verbatim}
\def\BibTeX{{\rm B\kern-.05em{\sc i\kern-.025em b}\kern-.08em
    T\kern-.1667em\lower.7ex\hbox{E}\kern-.125emX}}
\usepackage{balance}

\usepackage{xcolor}
\usepackage{booktabs}
\usepackage{multirow}
\usepackage{multicol}
\usepackage{siunitx}
\usepackage{enumitem}
\begin{document}
\title{Learning-based State Estimation in Distribution Systems with Limited Real-Time Measurements}

%% To specify the authors when (number of affiliations <= 2)
% \author{
% \IEEEauthorblockN{Author n.1 Name per Affiliation A\\ Author n.2 Name per Affiliation A}
% \IEEEauthorblockA{(Affiliation A) Department Name of Organization \\
% Name of the organization, acronyms acceptable\\
% City, Country\\
% \{email author n.1, email author n.2\}@domain (if desired)}
% \and
% \IEEEauthorblockN{Author n.1 Name per Affiliation B\\ Author n.2 Name per Affiliation B}
% \IEEEauthorblockA{(Affiliation B) Department Name of Organization \\
% Name of the organization, acronyms acceptable\\
% City, Country\\
% \{email author n.1, email author n.2\}@domain (if desired)}
% }

%% To specify the authors when (number of affiliations > 2)
% \author{\IEEEauthorblockN{J. G. De la Varga\IEEEauthorrefmark{1},
% S. Pineda\IEEEauthorrefmark{1},
% J. M. Morales\IEEEauthorrefmark{1} and
% Á. Porras\IEEEauthorrefmark{1}}
% \IEEEauthorblockA{\IEEEauthorrefmark{1}Oasys Research Group, University of Málaga\\ 
% \{josegv,spineda,juan.morales,alvaroporras\}@uma.es}
% \IEEEauthorblockA{\IEEEauthorrefmark{2} Department Name of Organization B\\
% Name of the organization B,
% Address B\\ Emails if wanted}
% \IEEEauthorblockA{\IEEEauthorrefmark{3} Department Name of Organization C\\
% Name of the organization C,
% Address C\\ Emails if wanted}
% \IEEEauthorblockA{\IEEEauthorrefmark{4}Department Name of Organization D\\
% Name of the organization D,
% Address D\\ Emails if wanted}
% }
\author{J. G. De la Varga, \IEEEmembership{Student Member, IEEE}, 
S. Pineda, \IEEEmembership{Senior Member, IEEE}, J. M. Morales \IEEEmembership{Senior Member, IEEE}, and \'A. Porras
\thanks{This work was supported by the Spanish Ministry of Science and Innovation (AEI/10.13039/501100011033) through project PID2020-115460GB-I00. The work of J.~G.~De la Varga was supported by the Spanish Ministry of Science and Innovation training program for PhDs with fellowship number PRE2021-098958. \emph{Corresponding author: Salvador Pineda.}}
\thanks{The authors are with Oasys Research Group, University of M\'alaga Spain, 29071 M\'alaga, Spain (e-mail: spineda@uma.es).}
}

\maketitle

\vspace{-20pt}
\begin{abstract}
The task of state estimation in active distribution systems faces a major challenge due to the integration of different measurements with multiple reporting rates. As a result, distribution systems are essentially unobservable in real time, indicating the existence of multiple states that result in identical values for the available measurements. Certain existing approaches utilize historical data to infer the relationship between real-time available measurements and the state. 
Other learning-based methods aim to estimate the measurements acquired with a delay, generating pseudo-measurements. Our paper presents a methodology that utilizes the outcome of an unobservable state estimator to exploit information on the joint probability distribution between real-time available measurements and delayed ones. Through numerical simulations conducted on a realistic distribution grid with insufficient real-time measurements, the proposed procedure showcases superior performance compared to existing state forecasting approaches and those relying on inferred pseudo-measurements.  
\end{abstract}


\vspace{-0.5cm}
\section{Introduction}\label{sec:intro}
\IEEEPARstart{T}{he} task of state estimation plays a critical role in power systems as it enables accurate monitoring and assessment of the system's operating conditions. However, in distribution systems, state estimation faces a major challenge due to the inherent lack of real-time observability, as certain measurements can only be acquired with a delay \cite{gomez-exposito2015}. Furthermore, an electric power system may lose part of its input data due to topology changes, sensor failures or communication errors \cite{abur2004}, malicious attacks \cite{gao2022}, and electrical blackouts \cite{wang2016}. As a result, distribution systems are practically unobservable in real time, indicating the existence of multiple states that result in identical values for the available measurements. This effect is even aggravated because distribution systems are historically measurement-scarce and have to deal with multiple reporting rates and possible asynchronization \cite{cheng2023survey}. Meanwhile, the progress in distributed energy resources (DERs) and the necessity for actively managing their  operational tasks have resulted in the emergence of active distribution networks (ADNs) \cite{Papadopoulos22Distributed}. Central to this paradigm shift is the imperative for real-time monitoring, analysis, and control of grid operations, for which the state estimation is a key component.

Under these circumstances, supervisory control and data acquisition (SCADA) systems and feeder terminal units (FTUs) play an essential role in providing real-time measurements to the distribution management system. Additionally, micro-phasor measurement units ($\mu$-PMUs) \cite{liu2020d} and advanced metering infrastructure \cite{mohassel2014survey}, which include smart meters (SMs) \cite{wang2018review}, also have great potential to provide data for state estimation. Given the heterogeneous nature and notably slow reporting times of this measurement set, most existing works discard dynamic state estimation and favor a static state estimation framework over a forecasting-aided one \cite{cheng2023survey}. 

Hence, to obtain valid estimations in a real-time fashion, a static state estimation process must incorporate additional information as available measurements are usually insufficient to deliver observability. Many existing approaches aim at generating the pseudo-measurements required to make the system observable by leveraging known information \cite{dehghanpour2019}, either employing probabilistic and statistical methods or learning-based techniques. 
% Within the latter, machine learning techniques such as regression models or neural networks can be trained to learn the relationship between the available measurements and the delayed ones. These models can then be used to generate pseudo-measurements for the estimation process whereby the state variables of a system are inferred. 
One example of the former group is \cite{singh2010}, where the first approach to incorporate loads as pseudo-measurements in distribution system state estimation is introduced. 
% Gaussian mixture models (GMM) are used to capture temporal load correlation at the expense of being sensitive to parameter selection and computationally expensive for high-dimensional learning. 
Another example of statistical methods to generate pseudo-measurements is proposed in  \cite{zhao2018}, where spatial and temporal correlations are accounted for by means of a first-order vector autoregressive model. The latter group of learning-based techniques has aroused great interest in academia and industry, mainly through advances in Deep Learning and Neural Networks (NN) models. In \cite{manitsas2012}, an artificial NN is used to generate power injection pseudo-measurements from a few real measurements. 
%As the resulting error is non-Gaussian and thus, not suitable for weighted least squares state estimation, the associated error is decomposed into various components through GMM. 

However, previous models have some limitations, namely, the need for abundant high-quality training data, the black-box nature, and the risk of physically infeasible results \cite{huang2023}. These drawbacks could be partly averted by linking physics-based prior knowledge to the models. A physics-conditioned optimization problem is formulated in \cite{kamal2022} to find the optimal choice for injection pseudo-measurements given a few real-time measures, while considering that the probability distribution for such unknown power injections follows a previously trained Generative Adversarial Network model. Another example can be found in \cite{wang2020}, where an autoencoder framework is physics-enhanced for pseudo-measurement generation. This hybrid learning method employs Deep NNs to learn the state correlations while considering physical flows in a power system.

Alternatively, in contrast to augmenting conventional methods with pseudo-measurements, some references propose directly forecasting the state using the available information. Authors in \cite{m.vinodkumar1999} simulate the behavior of power systems in two steps: first predicting the state based on previous estimations, followed by a filtering of the new set of measurements; both models, implemented by two distinct NN. From the same family but from a different category, a Deep NN is used in \cite{mestav2019} to learn the mappings between real-time measurements and the state, in combination with a Bayesian state estimation algorithm for improved bad-data detection. State forecasting can also be useful to estimate the state of observable systems. For instance, the authors of \cite{zamzam2019} train a shallow NN to learn how to initialize the Gauss-Newton algorithm by mapping available measurements to a point near the actual latent state, using historical or simulation data.
% A different approach divides a system without full observability into measured buses with varying injections and not measured buses with stable conventional loads \cite{bhela2018}. By leveraging these assumptions, a novel method is introduced to solve the nonlinear power flow equations by coupling them over consecutive time periods.

Likewise,  integrating the physics into the models can also be used to improve state forecasting. In \cite{pagnier2021}, a graph neural network is used to reconstruct physical parameters such as admittances of power lines without the need for complete observability. Once the model is tuned, the state can be directly predicted through available measurements. Another NN is introduced in \cite{zhang2019}, whose physics-informed approach is included by mirroring the structure of a regularized and ‘locally linearized’ state estimation algorithm, showcasing better performance and efficient training. Lastly,  \cite{dabush2023} proposes an original method that assumes that the grid state vector is a smooth graph signal with respect to the system admittance matrix. Subsequently, the authors develop a regularized weighted least squares state estimator that takes into account the graph structure of the network and does not require full observability.

% Tackling the problem from a different framework, dynamic state estimation addresses the challenge of observability by incorporating measurements from different time instants \cite{Zhao2019Power}. However, it requires fast scan rates to capture system dynamics and relies on restrictive assumptions about the system's stationary nature.

% When the system configuration allows for the construction of a measurement matrix of low-rank, sparse signal recovery methods offer another avenue by employing matrix completion techniques to estimate the state of the network \cite{Liu2019Matrix}. However, this assumption may not hold true in all system configurations, especially when considering spatial correlations between loads at neighboring buses.

Against this background, this paper introduces a novel data-driven static state estimation methodology specifically designed to deal with a measurement set with multiple reporting rates and insufficient real-time measurements. One of the key advantages of our proposal is its unique utilization of an unobservable state estimator's outcome to extract valuable information regarding the joint probability distribution between real-time available measurements and delayed measurements. This innovative approach combines the benefits of physics-informed methods while offering a distinct advantage over existing approaches by providing a simple, transparent, and replicable solution. To demonstrate the effectiveness of our proposal, we conduct numerical simulations on a realistic case study with limited real-time measurements. The results of these simulations showcase the improvement of our proposed methodology when compared to existing state forecasting approaches and those that rely on inferred pseudo-measurements. 

The rest of the paper is organized as follows. In Section \ref{sec:methodology} we present the state estimation formulation and its particularities when it comes to unobservable systems. Afterward, our proposal and its advantages are described. In Section \ref{sec:numerical results} we provide the numerical results of a case study to compare the performance of the proposed methodology with existing ones. Finally, Section \ref{sec:conclusion} provides concluding remarks.
% and outlines potential avenues for future research. 

\section{Methodology} \label{sec:methodology}
In this section, we begin by introducing the conventional static state estimation problem and outlining its key components. Subsequently, we delve into the challenges that arise when the observability of the entire network is not achieved, mainly focusing on state forecasting methods and pseudo-measurement generation techniques. Finally, we present our proposal that effectively improves upon existing methods by leveraging the information on the unobservable states provided by the system model equations and enhancing the accuracy of the overall estimation process. 

%First, we adopt a static state estimation framework and assume that the system parameters and measurements remain relatively constant over time. \textcolor{red}{Include something about dynamic state estimation here or in the introduction?}. 

Using a generic notation, the relationship between the measurement vector $z \in \mathbb{R}^m$ and the system state $x \in \mathbb{R}^{n}$ is characterized as:
\begin{equation}\label{eq:z_x_relation}
z=h(x)+e    
\end{equation}
\noindent where function $h(\cdot):\mathbb{R}^{n} \rightarrow \mathbb{R}^{m}$ is the measurement function that depends on the physical laws governing the relationship between the measured magnitudes and the state variables, and $e\in \mathbb{R}^{m}$ denotes the measurement error vector, whose components are typically assumed to follow normal distributions with zero mean and be independent and uncorrelated with respect to the state variables \cite{abur2004}. 

If the number of measurements is greater than the number of unknown state variables, and these measurements are sufficiently redundant and diverse, the system satisfies the observability requirements and the state estimation problem can be formulated as the following Weighted Least Square (WLS) optimization problem:
\begin{equation} 
\label{eq:wls-se} \hat {x} \in \arg \min_x
\left( z - h \left(x\right)\right)^{T}W\left( z - h \left(x\right)\right)
\end{equation}
\noindent where $\hat {x}$ is the estimated state vector, $(\cdot)^T$ is the matrix transposition operation, and the weight matrix $W=diag\{\sigma^{-2}_1,\ldots,\sigma^{-2}_m\}$) represents the user's confidence in the measured data, where $\sigma^2_m$ represents the error variance of the $m$-$th$ measurement. For ease of notation, the solution of problem \eqref{eq:wls-se} is denoted as $\hat{x}=\texttt{wls}(z)$. In the general case, problem \eqref{eq:wls-se} results in a nonconvex optimization problem without constraints, which conventionally has been solved using an iterative Newton-Raphson algorithm  \cite{wang2019}. Alternatively, problem \eqref{eq:wls-se} can also be directly solved by the interior-point methods implemented in nonlinear optimization solvers \cite{caro2020state}.

As discussed in Section \ref{sec:intro}, even when the numbers and types of measures captured by the sensors are enough to make the system fully observable, some of the measurements are not available in real time due to delays or communication issues. Consequently, we denote the set of available real-time measures as $z^a\in\mathbb{R}^{m_a}$, and the set of delayed measures as $z^d\in\mathbb{R}^{m_d}$, and split the measurement function as follows:
\begin{subequations}\label{eq:za_zd_x_relation}
\begin{align}
& z^a=h^a(x)+e^a \\
& z^d=h^d(x)+e^d
\end{align} 
\end{subequations}
where $h^a(\cdot):\mathbb{R}^{n} \rightarrow \mathbb{R}^{m_a}$, $h^d(\cdot):\mathbb{R}^{n} \rightarrow \mathbb{R}^{m_d}$, $e^a\in \mathbb{R}^{m_a}$, $e^d\in \mathbb{R}^{m_d}$, and $m=m_a+m_d$. In this work, we assume that the number and type of measures available in real time are not enough to make the system observable \cite{gomez-exposito2015, cheng2023survey}. Furthermore, we consider we have access to a collection of historical data comprising both available and delayed measurements. We denote this data set as $\{(z^a_k,z^d_k)\; \forall k \in \mathcal{K} \}$. Although model \eqref{eq:wls-se} cannot be used to infer the system state in real time because it is not observable, we can always perform a retrospective state estimation using all the historical measurements and data collected. That is, for each instance $k$ the retrospective system state denoted as $\hat{x}^r_k$ can be computed as $\hat{x}^r_k=\texttt{wls}(z^a_k,z^d_k)$. Therefore, we assume we have access to the enlarged historical data set $\{(z^a_k,z^d_k,\hat{x}^r_k)\; \forall k \in \mathcal{K} \}$. Under this framework, the main goal is to find a $\omega$-parameterized function $f(\cdot; \omega)$ that provides, as accurately as possible, the retrospective system state $\hat{x}^r_k$ out of the measurements available in real time $z^a_k$. Formally, we can estimate this function as follows:
\begin{equation}\label{eq:optimal_paramaters}
\omega^* \in \arg \min_\omega \sum_{k\in\mathcal{K}} (\hat{x}^r_k - f(z^a_k;\omega))^2
\end{equation}
Once parameters $\omega$ are properly tuned, the system state of a new instance $k'$ can be directly determined as
\begin{equation}\label{eq:direct}
\hat{x}_{k'} = f(z^a_{k'};\omega^*)
\end{equation}

Depending on the characteristics of function $f$, different strategies can be derived. This paper examines existing strategies that fall within this general framework while introducing novel approaches that surpass them in terms of estimation accuracy. In the simplest case, we can use the historical data to infer the joint probability distribution of the available measurements $z^a_k$ and the retrospective states $\hat{x}^r_k$ using conventional supervised regression techniques. These approaches are usually known as \emph{state forecasting} and multiple examples can be found in the literature, particularly since the tremendous boom of NN for prediction tasks \cite{m.vinodkumar1999, zamzam2019}.

In most cases, state forecasting methods consider standard functions $f$ that either simplify the underlying physical relationship between state variables and measurements in order to obtain interpretable models, or try to approximate function $h$ with elaborated machine learning approaches. For instance, if $f$ belongs to the family of linear functions, model \eqref{eq:optimal_paramaters} reduces to a linear regression problem that can be solved analytically to determine the optimal parameters $\omega^*$. However, the assumption of a linear relationship between state variables and measurements is highly restrictive and can significantly compromise the accuracy of this approach. On the other end of the spectrum, function $f$ can be represented by a neural network, providing a more flexible and adaptable modeling approach at the expense of losing interpretability, increasing sensitivity to parameter tuning, or augmenting the risk of overfitting. Finally, some recent research works propose physics-informed learning models to forecast the system state taking into account the underlying electricity network. However, the promising results in \cite{mestav2019} are derived from the careful tuning of some of the design parameters for the specific systems studied, while the method in \cite{zhang2019} exploits a NN with an embedded ‘locally linearized’ representation of the state-estimation equations, with the consequent approximation error.

Alternatively, other methodologies aim at inferring the joint probability distribution between available and delayed measurements using a function $g(\cdot; \lambda)$, parameterized on $\lambda$. For the same historical data set $\{(z^a_k,z^d_k,\hat{x}^r_k)\; \forall k \in \mathcal{K} \}$, the optimal value of $\lambda$ can be determined as follows:
\begin{equation}\label{eq:optimal_paramaters_pseudo}
\lambda^* \in \arg \min_\lambda \sum_{k\in\mathcal{K}} (z^d_k - g(z^a_k;\lambda))^2
\end{equation}
For an unseen instance $k'$, the delayed measurements can be estimated in real time as $\hat{z}^d_{k'} = g(z^a_{k'};\lambda^*)$. Since $\hat{z}^d_{k'}$ are not actual measurements provided by sensors, they are known as pseudo-measurements in the technical literature \cite{primadianto2016review}. By leveraging techniques for generating pseudo-measurements, the system becomes observable and its state can be determined as $\hat{x}_{k'}=\texttt{wls}(z^a_{k'},\hat{z}^d_{k'})$. Therefore, using a common mathematical framework, the relation between the system state and the available measurements is characterized as
\begin{equation} \label{eq:pseudo}
\hat{x}_{k'} = \texttt{wls}(z^a_{k'}, g(z^a_{k'};\lambda^*))   
\end{equation}
Again, selecting the proper function $g$ and solving problem \eqref{eq:optimal_paramaters_pseudo} are the main drawbacks of this strategy. Moreover, akin to the state forecasting approach discussed earlier, most existing methodologies generate pseudo-measurements either with a simple physical-agnostic approach or with complex physics-informed strategies as discussed in Section \ref{sec:intro}.

In this paper, we present a methodology aimed at improving the performance of state forecasting approaches and pseudo-measurement generation techniques in cases where the available real-time measurements $z^a$ are insufficient to ensure system observability. Specifically, we seamlessly integrate the measurement function $h$ into these techniques to enhance the estimation accuracy of the system state variables or the pseudo-measurements, respectively. To do so, for each instance $k$, we formulate the following optimization problem 
\begin{equation} 
\label{eq:wls_no_observable} \tilde{x}_k \in \arg \min_x
\left( z^a_k - h^a \left(x\right)\right)^{T}W^a\left( z^a_k - h^a \left(x\right)\right)
\end{equation}
where $W^a$ denotes the weight matrix of available measurements only. Since the number of available measurements is lower than the number of state variables, problem \eqref{eq:wls_no_observable} has multiple solutions and its optimal value is always 0. The optimization solver, therefore, delivers one solution satisfying $z^a_k = h^a(\tilde{x}_k)$, and thus, \emph{carrying} critical information on the measurement function $h$. Once problem \eqref{eq:wls_no_observable} is solved, we can generate pseudo-measurements as $\tilde{z}^d_k = h^d(\tilde{x}_k)$. This procedure guarantees that there exists a state, namely, $\tilde{x}_k$ that conforms to \emph{all} the measurements (both physical and pseudo-). Or in other words, this set of measurements is consistent with $h$. Conversely, if the pseudo-measurements $\tilde{z}^d_k$ are generated using the physics-agnostics methods previously discussed, there is no guarantee of the existence of state variables $x_k$ that simultaneously satisfy  $z^a_k = h^a(x_k)$ and $g(z^a_k;\lambda^*) = h^d(x_k)$. Additionally, the proposed scheme seems straightforward and transparent in contrast to the intricate physics-informed methods previously analyzed \cite{kamal2022} or \cite{zhang2019}. For notation purposes, we denote the solution of problem \eqref{eq:wls_no_observable} as $\tilde{x}=\texttt{un}(z^a)$.

The pseudo-measurements $\tilde{z}^d_k=h^d(\texttt{un}(z^a_{i}))$ we propose can be used in different ways. For instance, we can use them as additional features or explanatory variables in the state forecasting procedure formulated in \eqref{eq:optimal_paramaters} and \eqref{eq:direct}. By doing so, the optimal parameters $\omega$ are obtained as
\begin{equation}\label{eq:optimal_paramaters_improved}
\omega^* \in \arg \min_\omega \sum_{k\in\mathcal{K}} (\hat{x}^r_k - f(z^a_k,h^d(\texttt{un}(z^a_{k}));\omega))^2
\end{equation}
and the system state of a new instance $k'$ is computed as
\begin{equation}\label{eq:direct_improved}
\hat{x}_{k'} = f(z^a_{k'},h^d(\texttt{un}(z^a_{k'}));\omega^*)
\end{equation}

Note that the only difference between approaches \eqref{eq:optimal_paramaters}-\eqref{eq:direct} and \eqref{eq:optimal_paramaters_improved}-\eqref{eq:direct_improved} is that the proposed one also considers as explanatory variables or features of the regression model the pseudo-measurements $\tilde{z}^d_k$. In that sense, despite its apparent complexity, problem \eqref{eq:optimal_paramaters_improved} is a least-squared method with extra features obtained by solving problem \eqref{eq:wls_no_observable} for each instance $k$ with a non-convex optimization solver, and evaluating the known function $h^d$ on the obtained system state. However, unlike approach \eqref{eq:optimal_paramaters}-\eqref{eq:direct}, the estimation model \eqref{eq:optimal_paramaters_improved}-\eqref{eq:direct_improved} accounts for the measurement function $h$ through problem $\texttt{un}(\cdot)$ and through function $h^d$. Besides, this enhanced state forecasting approach is simple, direct and easy to interpret unlike the physics-informed black-box methodologies proposed in the literature.

Alternatively, the pseudo-measurements $\tilde{z}^d_k$  can also be used as features to learn the delayed measurements $z^d_k$ as follows
\begin{equation}\label{eq:optimal_paramaters_pseudo_improved}
\lambda^* \in \arg \min_\lambda \sum_{k\in\mathcal{K}} (z^d_k - g(z^a_k,h^d(\texttt{un}(z^a_{k}));\lambda))^2
\end{equation}

Similarly to \eqref{eq:optimal_paramaters_improved}, problem \eqref{eq:optimal_paramaters_pseudo_improved} is a conventional least-squared optimization model that can be easily solved once the proposed pseudo-measurements are computed. However, unlike problem \eqref{eq:optimal_paramaters_pseudo}, the one we propose incorporates physical information about the underlying electricity network. For a new instance $k'$, the pseudo-measurements required to make the system observable are obtained as
\begin{equation}\label{eq:pseudo_im}
\hat{z}^d_{k'} = g(z^a_{k'},h^d(\texttt{un}(z^a_{k'}));\lambda^*)   
\end{equation}

While the results obtained from this enhanced pseudo-measurement generation procedure can serve various purposes, their most intuitive application is as measurements in the conventional state estimation problem \eqref{eq:wls-se}. Consequently, the state variables can be determined as
\begin{equation} \label{eq:pseudo_improved}
\hat{x}_{k'} = \texttt{wls}(z^a_{k'},g(z^a_{k'},h^d(\texttt{un}(z^a_{k'}));\lambda^*))
\end{equation}
where the error variance $\sigma$ of the pseudo-measurements are obtained through the residuals of the regression model \eqref{eq:optimal_paramaters_pseudo_improved}.

In summary, this section describes four models, namely, \eqref{eq:direct}, \eqref{eq:pseudo}, \eqref{eq:direct_improved}, \eqref{eq:pseudo_improved}, which basically involve mapping functions between the available measurements $z^a$ and the state variables $x$. On the one hand, approaches \eqref{eq:direct} and \eqref{eq:direct_improved} represent state forecasting models that directly determine the system state using the available measurements. However, while \eqref{eq:direct} relies solely on observed data processed through physics-agnostic methods, the enhanced proposed model \eqref{eq:direct_improved} incorporates the system's physical equations through the measurement function $h$. On the other hand, approaches \eqref{eq:pseudo} and \eqref{eq:pseudo_improved} aim at generating pseudo-measurements to ensure system observability and determine the system state using standard estimation techniques. Correspondingly, while model \eqref{eq:pseudo} utilizes pseudo-measurements generated purely through physics-agnostic techniques, the proposed methodology \eqref{eq:pseudo_improved} does incorporate the physical behavior of the system. Furthermore, in contrast to alternative physics-informed methodologies prevalent in the technical literature, the models presented in this paper exhibit a straightforward, direct approach that facilitates effortless implementation and interpretation.

\section{Numerical Results} \label{sec:numerical results}
In the following, we provide results from a series of numerical experiments designed to test the performance of our state estimation approaches against alternative methods. We start by describing the setup of these experiments.

\subsection{Simulation setup}
\label{sec:setup}
\begin{figure}[ht!]
\centering
\includegraphics[width=3.3in]{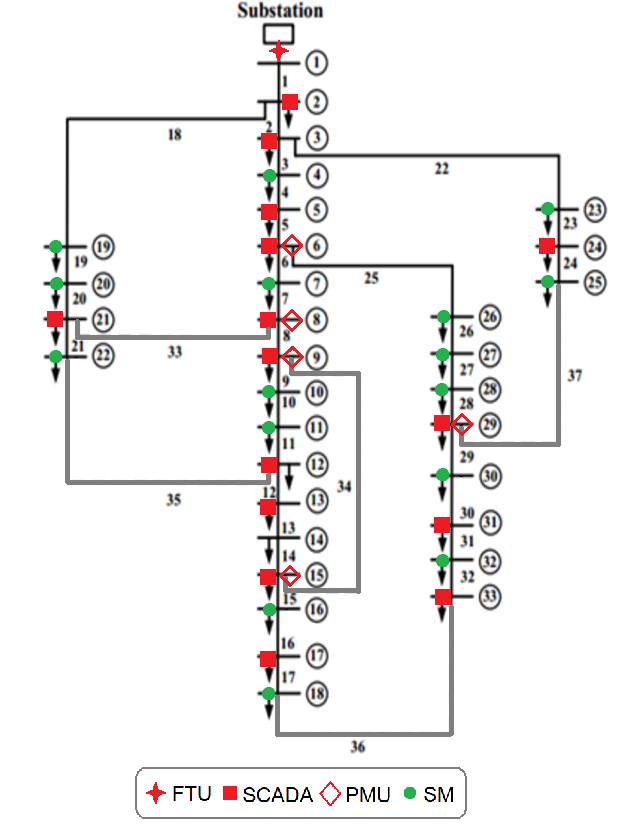}
\caption{IEEE-33 bus test case network and measurement set location.}
\label{fig:base_case_netplot}
\end{figure}

In the series of experiments, we utilize the well-known IEEE 33-bus distribution system test case \cite{Li23ieee33}. The data for this system is obtained from the Matpower test cases library \cite{zimmerman2010matpower}.
We have chosen the fully connected configuration of this system for a partially meshed topology, which consists of 33 nodes and 37 lines.
Fig.~\ref{fig:base_case_netplot} depicts the distribution network and the location of the measuring devices for some of the analysis discussed in this work.

To simulate a diverse set of operating conditions in ADNs with a high penetration of renewables, prosumers, storage, and other DERs \cite{Hidalgo10ADNreview}, the voltage magnitude ($V_i$) is randomly chosen within a range of 0.95 p.u. to 1.05 p.u. Similarly, the voltage phase angle ($\theta_i$) is selected between $-15\degree$ and $+15\degree$ \cite{massignan2022bayesian}. In state estimation, the measurement functions $h$ relate the set of measurable outputs to the system's state variables ($V_i,\theta_i$) as follows
\begin{subequations} \label{eq:h}
\begin{align} 
    &h_{V_i} = V_i \\
    &h_{\theta_i} = \theta_i \\    
    &h_{P_{ij}} = V_iV_j \left( G_{ij} \cos{\theta_{ij}} + B_{ij} \sin{\theta_{ij}} \right) - G_{ij}V_i^2 \\    
    &h_{Q_{ij}} = V_iV_j \left( G_{ij} \sin{\theta_{ij}} - B_{ij} \cos{\theta_{ij}} \right) + V_i^2\left(B_{ij}-b_{ij}^S/2\right)\\
    &h_{P_i} = V_i \sum _{j} V_j \left( G_{ij} \cos{\theta_{ij}} + B_{ij} \sin{\theta_{ij}} \right) \\
    &h_{Q_i} = V_i \sum _{j} V_j \left( G_{ij} \sin{\theta_{ij}} - B_{ij} \cos{\theta_{ij}} \right)
\end{align}
\end{subequations}
where $G_{ij}$ and $B_{ij}$ are, respectively, the real and imaginary parts of the bus admittance matrix, $b_{ij}^S$ is the shunt susceptance of the line $ij$, and ${\theta_{ij}=\theta_i - \theta_j}$. The system measurements are obtained by adding a Gaussian noise $\mathcal{N}(0,0.01^2)$ to power injections and flows and $\mathcal{N}(0,0.001^2)$ to voltage magnitudes and angles, as in \cite{zamzam2019}.

In state estimation for unobservable systems, measurement functions \eqref{eq:h} become particularly critical. For instance, while each voltage phasor measurement relates solely to an individual state variable, power flow measurements implicate the four state variables of the two buses of the line, and power injection measurements involve the state variables of all interconnected buses. Therefore, having access to measurements with functions that entail a high number of state variables offers richer information and provides more clues about the system's overall state, aiding in reconstructing the unobservable variables.

The measurement set illustrated in Fig.~\ref{fig:base_case_netplot} is detailed as follows.
\begin{itemize}
    \item Available measurements in real time ($m_a$): 
    \begin{itemize}
        \item One FTU provides  measurements of voltage magnitude and outgoing active and reactive power flows in the substation (red star).
        \item Fifteen active and reactive power injections are  available through SCADA measurements (red squares). 
        \item Five PMUs provide voltage magnitude and angle measurements (red diamonds). 
        %They are placed to be evenly distributed along the network in order to achieve better observability . 
    \end{itemize} 
    \item Delayed measurements ($m_d$): 
    \begin{itemize}
        \item Eighteen Smart Meters  placed in the remaining nodes provide power injection measurements (green circles). As previously stated, they are not available in real time due to their slower reporting times \cite{cheng2023survey}.
    \end{itemize}     
\end{itemize}

In real time we only have access to the set $m_a=3+2\times(15+5)=43$ while the measurements of $m_d=2\times(18)=36$ are delayed. Since $n=65$ \,(we fix $\theta_1$ to 0), the system is unobservable in real time and its state cannot be estimated through conventional techniques. 
However, the retrospective state estimation is determined assuming that we have both sets of available and delayed measurements, that is, $m=79$ measurements.

In Section \ref{sec:methodology} we describe two learning-based procedures to forecast the system state using function $f$ and other two approaches to generate pseudo-measurements through function $g$. Due to its simplicity and explainability, we compare in this case study the performance of these four approaches assuming that functions $f$ and $g$ belong to the family of linear functions. To learn these functions, we randomly generate 10\,000 instances. This set is then divided into two subsets: a training dataset comprising 8\,000 instances and a test dataset consisting of 2\,000 instances. By splitting the data in this manner, we can effectively train and evaluate the performance of the four approaches discussed.

The state estimation problems \eqref{eq:wls-se} and \eqref{eq:wls_no_observable} are modeled in Pyomo 6.4.2 \cite{pyomo} running in Python 3.10.13 and solved with Ipopt 3.12.8 \cite{ipopt}. In all cases, the initial solution of the algorithm is the flat state, that is, all voltage magnitudes and angles are set to 1 p.u. and $0\degree$, respectively. The adjustment of the parameters of functions $f$ and $g$ are determined using the linear regression function from the Scikit-learn package \cite{scikit-learn}. 

In order to demonstrate the performance of the proposed methodology, we compare in this section the computational results of the following nine approaches:
\begin{itemize}[leftmargin=*]
    \item[-] Benchmark (BN): This approach assumes that all measurements are available and the system is then observable. Under the considered setup, however, this situation is unrealistic and only used for benchmarking purposes.
    \item[-] Unobservable (UN): This approach consists of the solution of the optimization problem \eqref{eq:wls_no_observable} delivered by the non-convex optimization solver. Recall that, as previously discussed, problem \eqref{eq:wls_no_observable} has multiple solutions.
    \item[-] Flat (FL): This approach boils down to fixing the state variables to the flat voltage profile. 
    \item[-] State forecasting (SF): This strategy approximates the relationship between the available measures and the state variables using a linear function and therefore, implies the solution of \eqref{eq:optimal_paramaters} and \eqref{eq:direct}.
    \item[-] Enhanced state forecasting (SF$^*$): This approach is similar to the previous one but adds the additional enhanced features, i.e., $h^{d}(\texttt{un}(z^{a}))$, described in Section \ref{sec:methodology} and therefore, it implies the solution of \eqref{eq:optimal_paramaters_improved} and \eqref{eq:direct_improved}.
    \item[-] 1-NN: This is a naive approach that infers the system state as that corresponding to the nearest neighbor measured as the Euclidean distance between the available measurements. 
    \item[-] 20-NN: This is another naive forecasting approach where the system state is computed as the average of the state variables of the 20 nearest neighbors measured as the Euclidean distance between the available measurements. 
    \item[-] Pseudo-measurements (PM): This approach first generates pseudo-measurements through \eqref{eq:optimal_paramaters_pseudo} and \eqref{eq:pseudo} and then solves the state estimation problem \eqref{eq:wls-se}.
    \item[-] Enhanced pseudo-measurements (PM$^*$): This approach also generates pseudo-measurements through \eqref{eq:optimal_paramaters_pseudo_improved} and \eqref{eq:pseudo_im}, using the same additional features as SF$^*$. Then it obtains the system state with \eqref{eq:pseudo_improved}.     
\end{itemize}

\subsection{Enhancing effect of the proposed methodology}
\label{sec:enhancing_effect}

\begin{figure}
\centering
\includegraphics[width=3.3in]{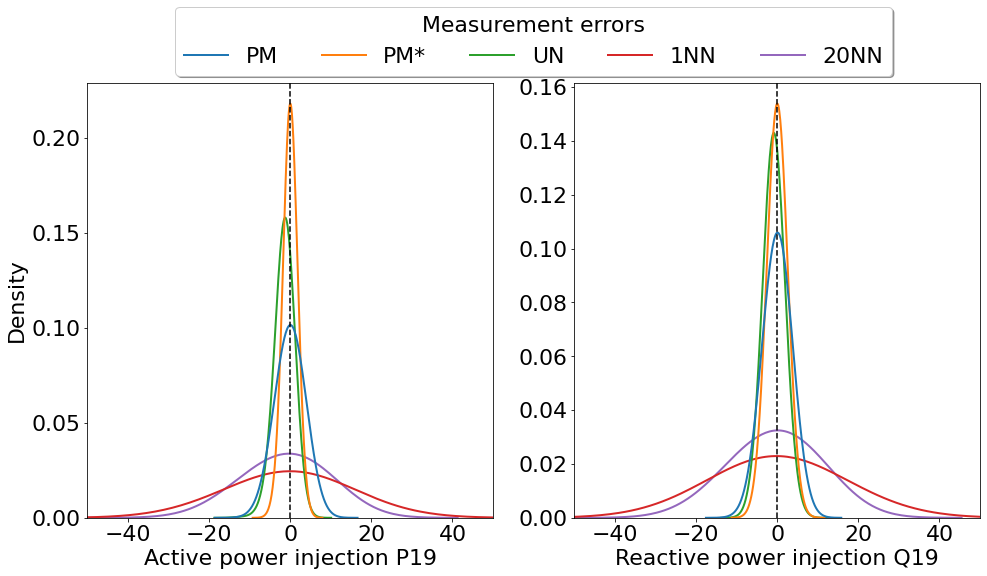}
\caption{Measurement residual distributions of power injections at node 19.}
\label{fig:dist_psm}
\end{figure}

Before presenting the computational results of the series of experiments, we illustrate in Fig.~\ref{fig:dist_psm} the enhancing effect of the pseudo-measurements $\tilde{z}^d$ obtained by solving the unobservable estimation problem \eqref{eq:wls_no_observable}. In particular, we plot the measurement error of the active and reactive power injection at node 19 for methods UN, PM, PM$^*$, 1-NN, and 20-NN in the test dataset, with respect to the measurements $z^d$. For the sake of illustration, the density plots have been depicted using a Gaussian kernel density estimate with a factor of 5 for the standard deviation of the smoothing kernel. As observed, although model \eqref{eq:wls_no_observable} has multiple solutions, the pseudo-measurements derived from the solution provided by the solver incur a lower variance in measurement residuals than some naive learning techniques such as 1-NN and 20-NN and in this node, even for the method PM. Additionally, using the yielded pseudo-measurements as features improves the performance of the regression model that aims at learning the delayed measures as a function of the available ones. For both the active and reactive power, we can observe that the residual probability distribution provided by PM$^*$ has a lower variance than that obtained by PM. 

The enhancing effect of the proposal is also shown in Table~\ref{tab:stats_pseudo}, which contains the mean and the standard deviation of the residuals of all delayed active power injections. The results in this table indicate that the pseudo-measurements obtained through the UN approach exhibit biased residuals. However, using these pseudo-measurements as features in the regression model of the  PM$^*$ approach effectively rectifies this issue, yielding unbiased estimations with the lowest variance for all delayed measurements corresponding to active power injections. This behavior is consistently replicated in all pseudo-measurements, but for brevity, only active power is presented.

The reduction in variance obtained through method UN is obviously reliant on the available measurement set. Solely relying on voltage phasor measurements (magnitude and angle) provides limited information. These measurements only constrain the corresponding state variables (bus voltages) to match their measured values, leaving the remaining unobservable state variables undetermined. In contrast, utilizing SCADA measurements rich in state variable connections, like power flows and power injections, offers a significant advantage since these measurements inherently involve the state variables of multiple buses. Consequently, the state estimator leads to a more constrained and informative solution, where the estimated state is not only consistent with the provided measurements but also closer to the true system state. This directly impacts the performance of the proposed approaches, SF* and PM*, as evidenced by the numerical results presented in this section.

\begin{table}[ht]
\renewcommand{\arraystretch}{1.3}
\centering
\caption{Active injection residuals statistics for the test set}
\label{tab:stats_pseudo}
\begin{tabular}{ccccccc}
\hline
\multirow{2}{*}{bus index} & \multicolumn{3}{c}{Mean, $\mu$} & \multicolumn{3}{c}{Standard deviation, $\sigma$} \\
 &    UN &      PM &     PM$^*$ &    UN &    PM &   PM$^*$ \\
\hline
% 1  &  0.00 &  0 &  0 &  0.014 &  0.014 &  0.014 \\
4  &  0.42 &  0 &  0 &  3.135 &  3.165 &  2.444 \\
7  & -0.11 &  0 &  0 &  0.849 &  1.201 &  0.771 \\
10 & -4.07 &  0 &  0 &  5.357 &  4.598 &  3.715 \\
11 &  3.35 &  0 &  0 &  6.724 &  5.358 &  4.425 \\
14 & -0.15 &  0 &  0 &  1.738 &  1.743 &  1.435 \\
16 & -0.60 &  0 &  0 &  1.921 &  1.648 &  1.583 \\
18 & 0 &  0 &  0 &  2.308 &  2.137 &  2.016 \\
19 & -1.39 &  0 &  0 &  1.798 &  2.653 &  1.252 \\
20 &  0.23 &  0 &  0 &  2.340 &  2.261 &  2.124 \\
22 & -0.38 &  0 &  0 &  2.041 &  1.974 &  1.917 \\
23 &  0.23 &  0 &  0 &  2.930 &  2.618 &  2.366 \\
25 & -0.38 &  0 &  0 &  2.983 &  2.479 &  2.188 \\
26 & -1.41 &  0 &  0 &  4.288 &  5.309 &  4.086 \\
27 & -0.25 &  0 &  0 &  5.110 &  5.120 &  5.009 \\
28 & -0.70 &  0 &  0 &  2.791 &  2.713 &  2.640 \\
30 & -0.58 &  0 &  0 &  3.270 &  2.254 &  2.198 \\
32 & -0.05 &  0 &  0 &  2.590 &  2.605 &  2.262 \\
\hline
\end{tabular}
\end{table}

\subsection{Results and discussion}
\label{sec:results}

In order to compare the performance of the nine approaches presented in Section \ref{sec:setup}, we use the Root Mean Squared Error (RMSE) defined as
\begin{equation}\label{metrics}
    \text{RMSE} = \sqrt{\frac{1}{K\cdot R}\sum_{k=1}^{K}\sum_{r=1}^{R}(y_{k,r} - \hat{y}_{k,r})^2}
\end{equation}
where $K$ is the number of instances in the test dataset, $R$ is either the number of buses or lines, depending on the magnitude of study; $y_{k,r}$ represents the true value of the magnitude, and $\hat{y}_{k,r}$ is the estimated value by each approach. Table \ref{tab:rmse_base_case} collates the RMSE values for the state variables (voltage magnitudes in p.u. and angles in degrees) and the power variables (active and reactive power injections and flows in MW and MVAr). 
%The power variables of each method are calculated by evaluating \eqref{eq:h} for each respective state. 
By disregarding the unrealistic results of the benchmark approach, the minimum error per magnitude is identified in bold.

\begin{table}[h]
\renewcommand{\arraystretch}{1.3}
\centering
\caption{RMSE values for the test set}
\label{tab:rmse_base_case}
\begin{tabular}{lcccccc}
\hline
Method & \textbf{$V$} & \textbf{$\theta$} & \textbf{$P$} & \textbf{$Q$} & \textbf{$P_f$} & \textbf{$Q_f$}\\
\hline
BN	&	0.004	&	0.003	&	0.11	&	0.11	&	0.24	&	0.24	\\
FL  &   0.029   &   8.641   &   76.96   &   114.04  &   41.88   &   63.43   \\
UN	&	0.052	&	0.089	&	25.85	&	33.73	&	17.72	&	22.49	\\
SF	&	0.021	&	0.086	&	22.03	&	32.84	&	14.53	&	21.56	\\
SF*	&	\textbf{0.018}	&	\textbf{0.085}	&	20.06	&	31.77	&	\textbf{13.40}	&	\textbf{20.96}	\\
PM	&	0.085	&	0.097	&	21.79	&	32.76	&	15.87	&	22.24	\\
PM*	&	0.041	&	0.087	&	\textbf{19.94}	&	\textbf{31.59}	&	13.48	&	20.97	\\
1NN	&	0.038	&	0.170	&	63.55	&	99.33	&	38.01	&	59.36	\\
20NN	&	0.027	&	0.120	&	44.54	&	68.16	&	26.76	&	41.08	\\

% Sf in p.u.
% BN	&	0.004	&	0.003	&	0.01	&	0.01	&	0.02	&	0.02	\\
% FL  &   0.029   &   8.641   &   7.70    &   11.40   &   4.19    &   6.34    \\
% UN	&	0.052	&	0.089	&	2.59	&	3.37	&	1.77	&	2.25	\\
% SF	&	0.021	&	0.086	&	2.20	&	3.28	&	1.45	&	2.16	\\
% SF*	&	\textbf{0.018}	&	\textbf{0.085}	&	2.01	&	3.18	&	\textbf{1.34}	&	\textbf{2.10}	\\
% PM	&	0.085	&	0.097	&	2.18	&	3.28	&	1.59	&	2.22	\\
% PM*	&	0.041	&	0.087	&	\textbf{1.99}	&	\textbf{3.16}	&	1.35	&	\textbf{2.10}	\\
% 1NN	&	0.038	&	0.170	&	6.35	&	9.93	&	3.80	&	5.94	\\
% 20NN	&	0.027	&	0.120	&	4.45	&	6.82	&	2.68	&	4.11	\\
\hline
\end{tabular}
\end{table}

We begin our analysis of the results shown in Table \ref{tab:rmse_base_case} by considering the outcomes obtained through the UN method first. Although this approach estimates the state of an unobservable system, it significantly reduces estimation errors in power variables and voltage angles, surpassing the performance of the naive approaches 1-NN and 20-NN. Conversely, the voltage magnitude errors are increased with respect to these methods.
This conclusion is aligned with the results provided in Fig.~\ref{fig:dist_psm} and Table~\ref{tab:stats_pseudo} and clearly indicates that the state yielded by the non-unique solution problem \eqref{eq:wls_no_observable} contains valuable information to be exploited. 
In fact, the results in Table \ref{tab:rmse_base_case} also show the improvement achieved by incorporating the pseudo-measurements derived from the UN method as additional features in the state estimation procedures since both SF$^*$ and PM$^*$ exhibit the lower estimation errors in all variables than their respective standard versions, SF and PM. Furthermore, it is worth noting that SF$^*$ and PM$^*$ exhibit the lowest estimation errors for voltage and power injection variables, respectively. This makes sense as the regression problems tackled by SF$^*$ and PM$^*$ are designed to respectively minimize the squared error of the voltage and of a measurement set with high proportion of power injection measurements. This reveals that incorporating a solution of an unobservable state estimation routine as explanatory features into state forecasting or pseudo-measurement generation procedures has a 
% substantial 
positive effect on the state estimation accuracy.

\begin{figure}[h]
\centering
\includegraphics[width=3.3in]{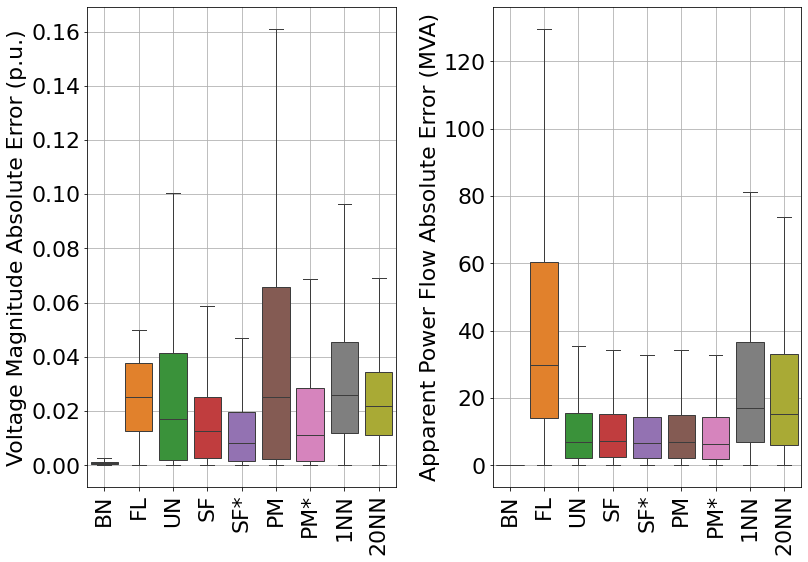}
\caption{Comparative of the absolute estimation error for the nine presented approaches of the voltage magnitude (left) and the apparent power flow (right).}
\label{fig:boxplot_base_case}
\end{figure}

Lastly, in order to provide a comprehensive view of the estimation errors, Fig.~\ref{fig:boxplot_base_case} complements Table~\ref{tab:rmse_base_case} with box plots corresponding to the distributions of the estimation errors for the nine compared approaches, with ${K\cdot R}$ data points in each case. Specifically, the voltage magnitude and the apparent power flow have been chosen as the key magnitudes. The voltage magnitude is crucial to ensure a stable operation and prevent equipment damage, while the precise estimation of apparent power flow is necessary to ensure that the power system operates within its thermal and electrical limits, avoiding overloads and potential failures. For ease of illustration, outliers have been removed from the plots. The figure reveals that the enhancements introduced in SF$^*$ and PM$^*$ not only reduce the mean error but also decrease its variance, with respect to SF and PM, respectively. This observation highlights the effectiveness of the proposed enhancing methodology in achieving more consistent and accurate state estimations.

The second part of this analysis investigates how the proposed methodology handles variations in active distribution network operating conditions.  These variations, often influenced by the penetration of renewables, storage, and prosumers, can significantly impact performance.  To evaluate this effect, state variables are generated using different uniform probability distribution ranges, as detailed in Table \ref{tab:variability} for low, medium, and high variability scenarios.

\begin{table}[h]
\renewcommand{\arraystretch}{1.3}
\centering
\caption{Variability scenarios}
\label{tab:variability}
\begin{tabular}{ccc}
\hline
Variability & Voltage Magnitude ($V$) & Voltage angle ($\theta$)\\
\hline
Low  &   $\mathcal{U}(0.975, 1.025)$ & $\mathcal{U}(-7.5\degree, +7.5\degree)$ \\
Medium & $\mathcal{U}(0.95, 1.05)$ & $\mathcal{U}(-15.5\degree, +15.5\degree)$ \\
High  &  $\mathcal{U}(0.925, 1.075)$ & $\mathcal{U}(-22.5\degree, +22.5\degree)$ \\
\hline
\end{tabular}
\end{table}

%The true state conditions formerly presented in section \ref{sec:setup} belong to active distribution networks and relate to the medium variability indicated in Table~\ref{tab:improvement_comparative}. In comparison, passive networks with a reduced penetration of distributed energy resources have a lower variability in their operating conditions, represented in the study by the low range $V_i~{\sim}~\mathcal{U}(0.975, 1.025)$ and $\theta_i~{\sim}~\mathcal{U}(-7.5\degree, +7.5\degree)$. For the sake of completeness we also present a situation with a higher variability range $(0.925, 1.075),\,(-22.5\degree, +22.5\degree)$. 

% Subsequent comparisons exclude the state forecasting method due to its lack of consideration for the system's physics, rendering it less reliable than the pseudo-measurement technique.
%

Beyond the variability of operating conditions, the performance of the proposed approaches can also be influenced by the type and number of available measurements, as discussed in Section \ref{sec:enhancing_effect}. To comprehensively evaluate the robustness of our proposed methodology, we conduct additional experiments where we vary the types and quantities of measurements used in the estimation process. In particular, the number of PMU and SCADA measurement devices is changed.

To quantify the performance of the proposed approaches across different scenarios, we  employ $\Delta$RMSE defined as the difference between the RMSE of the standard method and the enhanced one, that is,
\begin{equation}\label{metric_improvement}
    \Delta \text{RMSE} = \text{RMSE} - \text{RMSE}^*
\end{equation}

Consequently, positive values of this metric indicate that the enhanced method exhibits a lower error compared to the standard counterpart. Table~\ref{tab:improvement_comparative} summarizes the performance comparison between  PM and PM$^*$. Similar trends and conclusions were observed for the SF method, which are omitted for brevity. This table includes results for various variability scenarios (low, medium, high) and different measurement setups (number of PMU and SCADA devices). The table shows the difference in Root Mean Squared Error (RMSE) for both voltage magnitudes ($V$) and apparent power flow ($S_f$). Interestingly, all the values in Table~\ref{tab:improvement_comparative} are positive, signifying that the proposed method PM$^*$ consistently outperform PM across all scenarios.

%The final two columns of Table \ref{tab:improvement_comparative} present the accuracy improvement metric for the aforementioned critical magnitudes in distribution systems. The first column identifies if the improvement is evaluated either in the pseudo-measurement generation method or in the state forecasting technique. The second column categorizes the variability of operating conditions into three ranges. The third and fourth columns correspond to the number of PMUs and SCADA measurements, respectively.

\begin{table}
\renewcommand{\arraystretch}{1.3}
\renewcommand{\tabcolsep}{4pt}
\centering
\caption{$\Delta$RMSE for varying  variability and measurement data}
\label{tab:improvement_comparative}
\begin{tabular}{cccccc}
\hline
& Variability &	\# PMU	&	\# SCADA 	&	$\Delta \text{RMSE}\,(V)$ 	&	$\Delta \text{RMSE}\,(Sf)$	\\
\hline
% \multirow{6}{*}{\rotatebox[origin=c]{90}{SF vs. SF$^*$}} &	Low	&	5	&	15	&	0.001	&	0.093	\\
% &	Medium 	&	5	&	15	&	0.003	&	0.694	\\
% &	High	&	5	&	15	&	0.006	&	1.628	\\
% &	Medium 	&	15	&	5	&	0.002	&	0.209	\\
% &	Medium 	&	5	&	5	&	0.001	&	0.892	\\
% &	Medium 	&	5	&	25	&	0.009	&	0.173	\\
% \hline
\multirow{5}{*}{\rotatebox[origin=c]{90}{PM vs. PM$^*$}} &	Low	&	5	&	15	&	0.015	&	0.065	\\
&	Medium 	&	5	&	15	&	0.044	&	0.729	\\
&	High	&	5	&	15	&	0.070	&	1.957	\\
&	Medium 	&	15	&	5	&	0.002	&	0.070	\\
%&	Medium 	&	5	&	5	&	0.075	&	0.166	\\
&	Medium 	&	5	&	25	&	0.016	&	0.452	\\

% Sf in p.u.
% \multirow{6}{*}{PM} & Low	&	5	&	15	&	0.015	&	0.006	\\
% & Medium 	&	5	&	15	&	0.044	&	0.073	\\
% & High	&	5	&	15	&	0.070	&	0.196	\\
% & Medium 	&	15	&	5	&	0.002	&	0.007	\\
% & Medium 	&	5	&	5	&	0.075	&	0.017	\\
% & Medium 	&	5	&	25	&	0.016	&	0.045	\\
% \hline
% \multirow{6}{*}{SF} & Low	&	5	&	15	&	0.001	&	0.009	\\
% & Medium 	&	5	&	15	&	0.003	&	0.069	\\
% & High	&	5	&	15	&	0.006	&	0.163	\\
% & Medium 	&	15	&	5	&	0.002	&	0.021	\\
% & Medium 	&	5	&	5	&	0.001	&	0.089	\\
% & Medium 	&	5	&	25	&	0.009	&	0.017	\\
\hline
\end{tabular}
\end{table}

To assess the impact of variability on performance, we can examine the first three rows. It is apparent that the enhanced methods consistently outperform the baselines, with this advantage increasing  as variability rises. This trend suggests that the additional features introduced by the UN method become more valuable as operating conditions become less predictable.  In simpler terms, when conditions are stable, the standard methods using linear relationships perform well. However, under high variability, the additional features of the proposed methods help them handle the increased complexity of the system, leading to more accurate state estimation.

%Examining the initial three rows of each method in Table~\ref{tab:improvement_comparative}, it's clear that, despite varying operational conditions, the enhanced versions consistently outperform the baseline ones, albeit with different degrees of improvement. In both methods, the improvement increases with the variability of the true state, as the extra features provided by the UN method prove to be more beneficial. The explanation relies on the degree of non-linearity of the power flow equations, depending on the operation point. When variability is low, the relationship between power and state variables is nearly linear, making the linear regression present in $f$ and $g$ effective and rendering the extra features of PM* and SF* less useful. However, when variability is high, the non-linearities in power flow equations complicate estimation via linear regressions, making the additional features more useful. For a visual representation illustrating of how the non-linear effects of power flow equations vary with state variables, refer to Fig. 1 in \cite{coffrin2013linearprogramming}. 

%The second control variable of the numerical experiments is the set of measurements, as it directly impacts the state estimation problem formulation. Previous results were computed using the measurement set defined in section \ref{sec:setup} and in Fig.~\ref{fig:base_case_netplot}. For the last part of the discussion, we analyze the impact of modifying the available measurement set, both in size and shape. %type and number.

To understand how measurement type influences performance improvement, we compare rows 2 and 4 of Table~\ref{tab:improvement_comparative}.  These rows represent scenarios with the same total number of measurements (43) but different ratios of voltage phasor (PMU) and power injection (SCADA) devices.  As evident, the scenario with more SCADA devices (row 2) exhibits a higher improvement compared to the one with more PMUs (row 4). This aligns with the proposed methods' dependence on the UN method.  The UN method thrives on rich measurement sets that connect multiple state variables across different buses, providing more informative data compared to relying solely on isolated voltage measurements from PMUs.

Finally, we investigate the impact of the number of SCADA measurements, considering their previously observed advantage over PMU measurements. Interestingly, our analysis reveals an important trend when comparing rows 2 and 5 of Table~\ref{tab:improvement_comparative}.  While both scenarios utilize SCADA devices, a medium amount (15, as in row 2) provides greater improvement for both voltage magnitudes and apparent power flow compared to a high amount (25, as in row 5). This can be explained by the role SCADA measurements play in the pseudo-measurement generation process within method PM.  A larger set of SCADA injections significantly improves the estimation accuracy of method PM itself, reducing the need for the additional features introduced by the UN method.

\section{Conclusions} \label{sec:conclusion}

Our paper presents a novel approach to address the challenge of state estimation in distribution systems. This method leverages the outputs of an unobservable state estimator to exploit information about the joint probability distribution of real-time and delayed measurements.  It enriches learning-based methodologies by incorporating these outputs as additional features, enabling either direct state variable forecasting or real-time pseudo-measurement generation. Numerical simulations on a realistic distribution network demonstrate the superiority of our proposed method across various scenarios.  In particular, it demonstrates a clear advantage in active distribution networks with highly variable operating conditions. Additionally, the performance improves when the unobservable state estimator utilizes measurements encompassing a high number of state variables.

\section{Acknowledgments}
The authors thank the Supercomputing and Bioinnovation Center (SCBI) of the University of M\'alaga for their provision of computational resources (the supercomputer Picasso) and technical support (www.scbi.uma.es/site).

\bibliographystyle{IEEEtran}
\bibliography{references}

% thak's all folks
\end{document}